\documentclass{article}
\usepackage{amssymb}
\usepackage{latexsym}

\usepackage{graphicx}

\title{On Frankl and F\"{u}redi's conjecture for $3$-uniform hypergraphs}

\author{ Qingsong Tang \thanks{College of Sciences, Northeastern University, Shenyang, 110819, China and Mathematics School, Institute of Jilin University, Changchun, 130012, China. Email: t\_qsong@sina.com.cn} \and Hao Peng \thanks{ Corresponding author. College of Mathematics, Hunan University, Changsha 410082, China. Email: hpeng@hnu.edu.cn. Supported in part by National Natural Science Foundation of China (No. 11201135)} \and Cailing Wang \thanks{School of Mathematics, Jilin University, Changchun 130012, China. Email: wangcl-jl@163.com} \and Yuejian Peng \thanks{ College of Mathematics, Hunan University, Changsha 410082, China. Email: ypeng1@163.com. Supported in part by National Natural Science Foundation of China (No. 11271116).} }

\date{}
\newtheorem{defi}{Definition}[section]
\newtheorem{theo}{Theorem}[section]

\newtheorem{remark}[theo]{Remark}

\newtheorem{lemma}[theo]{Lemma}

\newtheorem{coro}[theo]{Corollary}
\newtheorem{con}[theo]{Conjecture}

\newtheorem{fact}[theo]{Fact}

\newcommand{\qed}{\hspace*{\fill} \rule{7pt}{7pt}}

\begin{document}
\maketitle
\begin{abstract}The Lagrangian of a hypergraph has been  a useful tool in hypergraph extremal problems. In most applications, we need an upper bound for the Lagrangian of a hypergraph.
Frankl and F\"{u}redi in \cite{FF} conjectured that the $r$-graph with $m$ edges formed by taking the first $m$ sets in the colex ordering of ${\mathbb N}^{(r)}$ has the largest Lagrangian of all $r$-graphs with $m$ edges. In this paper, we give some partial results for this conjecture.
\end{abstract}


\section{Introduction and the main results}

In 1965,Motzkin and Straus \cite{MS} provided a new proof of Tur\'an's theorem based on a continuous characterization of the clique number of a graph using Lagrangians of  graphs. This new proof aroused interests in the study of Lagrangians of $r$-graphs. The Lagrangian of a hypergraph has been  a useful tool in hypergraph extremal problems.

For a set $V$ and a positive integer $r$ we denote by $V^{(r)}$ the family of all $r$-subsets of $V$. An $r$-uniform graph or $r$-graph $G$ consists of a set $V(G)$ of vertices and a set $E(G) \subseteq V(G) ^{(r)}$ of edges. An edge $e=\{a_1, a_2, \ldots, a_r\}$ will be simply denoted by $a_1a_2 \ldots a_r$. An $r$-graph $H$ is  a {\it subgraph} of an $r$-graph $G$, denoted by $H\subseteq G$ if $V(H)\subseteq V(G)$ and $E(H)\subseteq E(G)$. The complement of an $r$-graph $G$ is denoted by $G^c$. Let $K^{(r)}_t$ denote the complete $r$-graph on $t$ vertices, that is the $r$-graph on $t$ vertices containing all possible edges. A complete $r$-graph on $t$ vertices is also called a clique with order $t$. Let ${\mathbb N}$ be the set of all positive integers. For an integer $n \in {\mathbb N}$, we denote the set $\{1, 2, 3, \ldots, n\}$ by $[n]$. Let $[n]^{(r)}$  represent the  complete $r$-graph on the vertex set $[n]$. When $r=2$, an $r$-graph is a simple graph.  When $r\ge 3$,  an $r$-graph is often called a hypergraph.

\begin{defi}
For  an $r$-graph $G$ with the vertex set $[n]$,
edge set $E(G)$ and a vector $\vec{x}=(x_1,\ldots,x_n) \in R^n$,
define
$$\lambda (G,\vec{x})=\sum_{i_1i_2 \cdots i_r \in E(G)}x_{i_1}x_{i_2}\ldots x_{i_r}.$$
\end{defi}

\begin{defi}
Let $S=\{\vec{x}=(x_1,x_2,\ldots ,x_n): \sum_{i=1}^{n} x_i =1, x_i
\ge 0 {\rm \ for \ } i=1,2,\ldots , n \}$. The Lagrangian of
$G$, denoted by $\lambda (G)$, is defined as
 $$\lambda (G) = \max \{\lambda (G, \vec{x}): \vec{x} \in S \}.$$
The value $x_i$ is called the {\em weight} of the vertex $i$. We call $\vec{x}=(x_1, x_2, \ldots, x_n) \in R^n$ a legal weighting for $G$ if $\vec{x} \in S$. A vector $\vec{y}\in S$ is called an {\em optimal weighting} for $G$ if $\lambda (G, \vec{y})=\lambda(G)$.
\end{defi}

 The following fact is easily implied by the definition of the Lagrangian.

\begin{fact}\label{mono}
Let $G_1$, $G_2$ be $r$-graphs and $G_1\subseteq G_2$. Then $\lambda (G_1) \le \lambda (G_2).$
\end{fact}

In \cite{MS}, Motzkin and Straus provided the following simple expression for the Lagrangian of a 2-graph.

\begin{theo} (Motzkin and Straus \cite{MS}) \label{MStheo}
If $G$ is a 2-graph in which a largest clique has order $t$ then
$\lambda(G)=\lambda(K^{(2)}_t)={1 \over 2}(1 - {1 \over t})$.
\end{theo}

The obvious generalization of Motzkin and Straus' result to hypergraphs is false because there are many examples of hypergraphs that do not achieve their Lagrangian on any proper subhypergraph.Lagrangians of hypergraphs has been proved to be a useful tool in hypergraph extremal problems. %Applications of Lagrangian method can be found in \cite{FR84}, \cite{FF}, \cite{mubayi06}, \cite{sidorenko89} and \cite{keevash}.
In most applications, an upper bound is needed.
Frankl and F\"uredi \cite{FF} asked the following question. Given $r \ge 3$ and $m \in {\mathbb N}$ how large can the Lagrangian of an $r$-graph with $m$ edges be?
For distinct $A, B \in {\mathbb N}^{(r)}$ we say that $A$ is less than $B$ in the {\em colex ordering} if $max(A \triangle B) \in B$, where $A \triangle B=(A \setminus B)\cup (B \setminus A)$. For example we have $246 < 156$ in ${\mathbb N}$ since $max(\{2,4,6\} \triangle \{1,5,6\}) \in \{1,5,6\}$. In colex ordering, $123<124<134<234<125<135<235<145<245<345<126<136<236<146<246<346<156<256<356<456<127<\cdots .$
Note that the first $t \choose r$ $r$-tuples in the colex ordering of ${\mathbb N}^{(r)}$ are the edges of $[t]^{(r)}$.

The following conjecture of Frankl and F\"uredi (if it is true) proposes a  solution to the question mentioned above.

\begin{con} (Frankl and F\"uredi \cite{FF})\label{conjecture} The $r$-graph with $m$ edges formed by taking the first $m$ sets in the colex ordering of ${\mathbb N}^{(r)}$ has the largest Lagrangian of all $r$-graphs with  $m$ edges. In particular, the $r$-graph with $t \choose r$ edges and the largest Lagrangian is $[t]^{(r)}$.
\end{con}

This conjecture is true when $r=2$ by Theorem \ref{MStheo}. For the case $r=3$, Talbot in \cite{T} proved the following.

\begin{theo} (Talbot \cite{T}) \label{Tal} Let $m$ and $t$ be integers satisfying
${t-1 \choose 3} \le m \le {t-1 \choose 3} + {t-2 \choose 2} - (t-1).$
Then Conjecture \ref{conjecture} is true for $r=3$ and this value of $m$.
Conjecture \ref{conjecture} is also true for $r=3$ and $m= {t \choose 3}-1$ or $m={t \choose 3} -2$.
\end{theo}
For the case $r=3$, Tang, Peng, Zhang, and Zhao in  \cite{TPZZ2} proved the following.
%\begin{theo} (\cite{TPZZ} \cite{TPZZ2}), \label{Tang} Let $m$ and $t$ be integers.Then  Conjecture \ref{conjecture} is  true for $r=3$  and $m= {t \choose 3}-3$, $m={t \choose 3} -4$ or $m={t \choose 3} -5$. Conjecture \ref{conjecture} is also true for $r=3$ and $m\leq 56$.
%\end{theo}

\begin{theo} \cite{TPZZ2} \label{Tang2} Let $m$ and $t$ be positive integers satisfying ${t-1 \choose 3} \le m \le {t-1 \choose 3} + {t-2 \choose 2}-(t-4)$. Then Conjecture \ref{conjecture} is true for $r=3$ and this value of $m$.
\end{theo}

The truth of Frankl and F\"uredi's conjecture is not known in general for $r \ge 4$. In the case $r=3$, the case when ${t-1 \choose 3}+{t-2 \choose 2}-(t-5) \le m \le {t \choose 3}-3$ is still open in this conjecture. Before we state our main results below, let us introduce some further notations and terminologies.

Let $C_{r,m}$ denote the $r$-graph with $m$ edges formed by taking the first $m$ sets in the colex ordering of ${\mathbb N}^{(r)}$.
Denote $$\lambda_{m}^{r}=\max\{\lambda(G): G {\rm \ is \ an \ } r-{\rm graph\ with \ } m {\rm \ edges }\}.$$

\begin{defi}
An $r$-graph $G$  with $m$ edges is called an extremal $r$-graph if  $\lambda(G)=\lambda_{m}^{r}$.
\end{defi}

Frankl and F\"uredi's conjecture says that $\lambda_{m}^{r}=\lambda(C_{r,m})$. To verify the truth of Conjecture \ref{conjecture}, it is sufficient to show that $\lambda(G)\le \lambda(C_{r,m})$ holds  for every  extremal $r$-graph $G$ with $m$ edges.

\begin{defi}
An $r$-graph $G=([n],E)$ is {\it left-compressed} if $j_1j_2 \cdots j_r \in E$ implies $i_1i_2 \cdots i_r \in E$ provided $i_p \le j_p$ for every $p, 1\le p\le r$. Equivalently, an $r$-graph $G=([n],E)$ is {\it left-compressed} if $E_{j\setminus i}=\emptyset$ for any $1\le i<j\le n$.
\end{defi}

The following lemma implies that we only need to consider left-compressed extremal $r$-graphs to verify that $\lambda(G)\le \lambda(C_{r,m})$ holds  for every  extremal $r$-graph $G$ with $m$ edges.

\begin{lemma} \label{lemmaleftcompress}\cite{T} There exists a left-compressed $r$-graph $G$ with $m$ edges  such that $\lambda(G)=\lambda_m^r$.
\end{lemma}

To emphasize this, let us make a remark.
\begin{remark}
To verify Conjecture \ref{conjecture}, it is sufficient to show that for a left-compressed extremal $r$-graphs $G$ with $m$ edges, $\lambda(G)\leq\lambda(C_{r,m})$ holds.
\end{remark}

Applying the following result showed in \cite{T}, we can further reduce the classes of $3$-graphs to verify in order to verify Conjecture \ref{conjecture}.

\begin{lemma} (\cite{T}) \label{LemmaTal9}Let $m$ be  a positive integer. Let $G=([n], E)$ be a left-compressed extremal $3$-graph with $m$ edges. If $\vec{x}=(x_{1},x_{2},\ldots ,x_{n})$  is an optimal weighting for $G$ satisfying $x_1 \ge x_2 \ge \ldots \ge x_k >x_{k+1}=\ldots=x_{n}=0$. Then
$$|E|\geq {k-1 \choose 3}+{k-2 \choose 2}-(k-2).$$
\end{lemma}

\begin{remark}\label{reducetot}
Let $G$ be a left-compressed extremal
$3$-graph with $m$ edges.  Let $t$ be a positive integer such that ${t-1 \choose 3} \le m < {t \choose 3}$. To show $\lambda(G)\leq\lambda(C_{3,m})$, we can assume $G$ is  on  vertex set $[t]$.
\end{remark}

\noindent{\em Proof.} Let $G=([n],E)$ and $t$ satisfy the conditions in this remark.  Let $\vec{x}=(x_{1},x_{2},\ldots ,x_{n})$  be an optimal weighting for $G$ satisfying $x_1 \ge x_2 \ge \ldots \ge x_k >x_{k+1}=\ldots=x_{n}=0$. We claim that $k\leq t$. Otherwise $k\geq t+1$ and Lemma \ref{LemmaTal9} implies that
\begin{eqnarray*}
m=|E|&\geq& {k-1 \choose 3}+{k-2 \choose 2}-(k-2)\nonumber\\
&\geq&{t \choose 3}+{t-1 \choose 2}-(t-1)\nonumber\\
&\ge &{t \choose 3}
\end{eqnarray*}
which contradicts to the assumption that ${t-1 \choose 3}\le m< {t \choose 3}.$ So to show $\lambda(G)\leq\lambda(C_{3,m})$, we can assume $G$ is on $[t]$.  \qed

In this paper, we first give the following partial result.

\begin{theo} \label{Lemma0} Let $m$, $t$, $a$ and $i$ be positive integers satisfying  $m={t \choose 3}-a$,  $3\leq a\leq t-2$ and $i\ge 1$. Let $G=([t],E)$ be a left-compressed $3$-graph  with $m$ edges. If the triple with minimum colex ordering in  $G^{c}$ is $(t-2-i)(t-2)t$. Then $\lambda(G)\leq\lambda(C_{3,m})$.
\end{theo}

In \cite{T}, the following result was proved.

\begin{theo}\label{Tal5}\cite{T}
Let m, t and a satisfy $-(t - 2) \le a \le (t - 5)$ and
$$m = {t-1 \choose 3} + {t - 2 \choose 2}+ a.$$
Suppose $G$ is a left-compressed extremal 3-graph with m edges. Then $G$ and $C_{3,m}$ differ in at most $2(t - a - 2)$ edges, i.e.,
$|E(G) \Delta E(C_{3,m})| \le  2(t - a - 2)$.
\end{theo}

 We show that

\begin{theo} \label{mainresult} Let $m$ be any positive integer. Let $G$ be a left-compressed  extremal $3$-graph with $m$ edges satisfying $|E(G) \Delta E(C_{3,m})| \le  6$. Then $\lambda(G)\leq\lambda(C_{3,m})$.
\end{theo}

In the proof of Theorem \ref{mainresult}, we will prove several lemmas in Section \ref{proof}. These lemmas themselves provide partial results to Conjecture \ref{conjecture} as well.

Using Theorem \ref{mainresult}, we can prove Conjecture \ref{conjecture} holds for ${t \choose 3}-6 \le m\le {t \choose 3}-3$ when $r=3$.
\begin{coro}\label{coro1}
Let $m$ and $t$ be  positive integers satisfying ${t \choose 3}-6 \le m\le {t \choose 3}-3$. Let $G$ be a $3$-graph with $m$ edges, then $\lambda(G)\leq \lambda (C_{3,m}).$
\end{coro}

All proofs will be given in Section \ref{proof}.

\section{Useful Results}

We will impose one additional condition on any optimal weighting ${\vec x}=(x_1, x_2, \ldots, x_n)$ for an $r$-graph $G$:
\begin{eqnarray}
 &&|\{i : x_i > 0 \}|{\rm \ is \ minimal, i.e. \ if}  \ \vec y {\rm \ is \ a \ legal \ weighting \ for \ } G  {\rm \ satisfying }\nonumber \\
 &&|\{i : y_i > 0 \}| < |\{i : x_i > 0 \}|,  {\rm \  then \ } \lambda (G, {\vec y}) < \lambda(G) \label{conditionb}.
\end{eqnarray}

For an $r$-graph $G=(V,E)$ we denote the $(r-1)$-neighborhood of a vertex $i \in V$ by $E_i=\{A \in V^{(r-1)}: A \cup \{i\} \in E\}$. Similarly, we will denote the $(r-2)$-neighborhood of a pair of vertices $i,j \in V$ by $E_{ij}=\{B \in V^{(r-2)}: B \cup \{i,j\} \in E\}$. We denote the complement of $E_i$ by $E^c_i=\{A \in V^{(r-1)}: A \cup \{i\} \in V^{(r)} \backslash E\}$. Also, we  denote the complement of $E_{ij}$ by
$E^c_{ij}=\{B \in V^{(r-2)}: B \cup \{i,j\} \in V^{(r)} \backslash E\}$. Denote $$E_{i\setminus j}=E_i\cap E^c_j.$$

When the theory of Lagrange multipliers is applied to find the optimum of $\lambda(G, {\vec x})$, subject to $\sum_{i=1}^n x_i =1$, notice that $\lambda (E_i, {\vec x})$ corresponds to the partial derivative of  $\lambda(G, \vec x)$ with respect to $x_i$. The following lemma gives some necessary conditions of an optimal weighting for $G$.

\begin{lemma} (Frankl and R\"odl \cite{FR84}) \label{LemmaTal5} Let $G=([n], E)$ be an $r$-graph  and ${\vec x}=(x_1, x_2, \ldots, x_n)$ be an optimal  weighting for $G$ with $k$   positive weights $x_1, x_2, \ldots, x_k$  satisfying condition (\ref{conditionb}). Then for every $\{i, j\} \in [k]^{(2)}$, (a) $\lambda (E_i, {\vec x})=\lambda (E_j, \vec{x})=r\lambda(G)$, (b) there is an edge in $E$ containing both $i$ and $j$.
\end{lemma}

\begin{remark}\label{r1} (a) In Lemma \ref{LemmaTal5}, part(a) implies that
$$x_j\lambda(E_{ij}, {\vec x})+\lambda (E_{i\setminus j}, {\vec x})=x_i\lambda(E_{ij}, {\vec x})+\lambda (E_{j\setminus i}, {\vec x}).$$
In particular, if $G$ is left-compressed, then
$$(x_i-x_j)\lambda(E_{ij}, {\vec x})=\lambda (E_{i\setminus j}, {\vec x})$$
for any $i, j$ satisfying $1\le i<j\le k$ since $E_{j\setminus i}=\emptyset$.

(b) If  $G$ is left-compressed, then for any $i, j$ satisfying $1\le i<j\le k$,
\begin{equation}\label{enbhd}
x_i-x_j={\lambda (E_{i\setminus j}, {\vec x}) \over \lambda(E_{ij}, {\vec x})}
\end{equation}
holds.  If  $G$ is left-compressed and  $E_{i\setminus j}=\emptyset$ for $i, j$ satisfying $1\le i<j\le k$, then $x_i=x_j$.

(c) By (\ref{enbhd}), if  $G$ is left-compressed, then an optimal legal weighting  ${\vec x}=(x_1, x_2, \ldots, x_n)$ for $G$  must satisfy
\begin{equation}\label{conditiona}
x_1 \ge x_2 \ge \ldots \ge x_n \ge 0.
\end{equation}
\end{remark}

We also need the following lemma from \cite{T} in the proof of our main results.

\begin{lemma} (Talbot \cite{T}) \label{LemmaTal7} For  integers $m,  t, $ and $r$ satisfying ${t-1 \choose r} \le m \le {t-1 \choose r} + {t-2 \choose r-1}$, we have $\lambda(C_{r,m}) = \lambda([t-1]^{(r)})$. \end{lemma}

%Denote $$\lambda_{m}^{r}=\max\{\lambda(G): G {\rm \ is \ an \ } r-{\rm graph\ with \ } m {\rm \ edges }\}.$$

%\begin{lemma} (Talbot \cite{T}) \label{LemmaTal7} For  integers $m,  t, $ and $r$ satisfying ${t-1 \choose r} \le m \le {t-1 \choose r} + {t-2 \choose r-1}$, we have $\lambda(C_{r,m}) = \lambda([t-1]^{(r)})$. \end{lemma}
%In \cite{T}, the following result is proved.

%\begin{lemma} (Talbot \cite{T}) \label{LemmaTal9}
%Let $G=(V,E)$ be a $3$-graph with $m$ edges such that $\lambda(G)=\lambda_{m}$. Let $\vec{x}=(x_{1},x_{2},\ldots ,x_{k})$  be an optimal weighting for $G$ satisfying $x_1 \ge x_2 \ge \ldots \ge x_k >x_{k+1}=\ldots=x_{n}=0$. Then $$|E|\geq {k-1 \choose 3}+{k-2 \choose 2}-(k-2).$$ \end{lemma}
%\begin{lemma}(Peng,Tang and Zhao \cite{PTZ})\label{Lemmaptz} Let $G$ be a left-compressed $r$-graph on the vertex set $[t]$ containing the clique $[t-1]^{(r)}$. Let ${\vec x}=(x_1, x_2, \ldots, x_t)$ be an optimal weighting for $G$. Then
%\begin{equation}
% x_1\le x_{t-1}+x_t \le 2x_{t-1}.
%\end{equation}
%\end{lemma}
%\begin{lemma}(Tang,Peng,Zhang and Zhao \cite{TPZZ2})\label{Lemmatpzz1} Let $G=(V,E)$ be a left-compressed 3-graph on the vertex set $[t]$. Let ${\vec x}=(x_1, x_2, \ldots, x_t)$ be an optimal weighting for $G$. Assume that $i_{1}\leq i_{2}$ and $j_{1}\leq j_{2}$ are integers. Then $\lambda(E_{i_{1}j_{1}},\vec{x})\geq \lambda(E_{i_{2}j_{2}},\vec{x})$.
%\end{lemma}

\section{Proof of Main Results}\label{proof}

\subsection{Proof of Theorem \ref{Lemma0}}

\noindent{\em Proof.}Since $G$ is left-compressed, in view of Figure 1, then we have $a\geq 2i+1$.

 To show that $ \lambda(G)\le \lambda(C_{3,m})$, we will take an optimal weighting $\vec{x}$ for $G$, then we take a legal weighting, say $\vec{z}$ for $C_{3,m}$ by replacing  a few coordinators of $\vec{x}$ and  show that $\lambda(G, \vec{x})\leq \lambda(C_{3,m}, \vec{z})$. This would imply that $$\lambda(G)=\lambda(G, \vec{x})\le\lambda(C_{3,m}, \vec{z})\leq \lambda(C_{3,m}).$$

 Let us go into the details.
Let $\vec{x}=(x_{1},x_{2},\ldots ,x_{t})$  be an optimal weighting for $G$ satisfying $x_1 \ge x_2 \ge \ldots \ge x_t \ge 0$. First we point out that
\begin{eqnarray}\label{ineqadd}
\lambda(E_{1(t-2-i)},\vec{x})-\lambda(E_{(t-2)(t-1)},\vec{x})&=& x_{t-2}+x_{t-1}+x_{t}-x_{1}-x_{t-2-i} \geq 0.
\end{eqnarray}
To verify (\ref{ineqadd}), by Remark \ref{r1}(b), we have
\begin{eqnarray}\label{By1}
x_{1}=x_{t-1}+\frac{\lambda(E_{1\backslash (t-1)},\vec{x})}{\lambda(E_{1(t-1)},\vec{x})}\leq x_{t-1}+\frac{(x_{2}+\cdots+x_{t-2})x_{t}}{x_{2}+\cdots+x_{t-2}+x_{t}}\leq x_{t-1}+x_{t};
\end{eqnarray}
 \begin{eqnarray}\label{By2}x_{1}&=&x_{t-2}+\frac{\lambda(E_{1\backslash (t-2)},\vec{x})}{\lambda(E_{1(t-2)},\vec{x})}\nonumber \\ &=&x_{t-2}+\frac{x_{t-2-i}+\cdots+ x_{t-3}+x_{t-1}}{1-x_{1}-x_{t-2}}x_{t}\nonumber \\
&\leq& x_{t-2}+\frac{x_{t-2-i}+\cdots+ x_{t-3}+x_{t-1}}{1-x_{t-2}-x_{t-1}-x_t}x_{t}  \quad \quad ({\rm By} \ (\ref{By1})) \nonumber \\
&\leq& x_{t-2}+\frac{x_{1}+x_2+\cdots+ x_{i}+x_{t-2}}{1-x_{t-2-i}-x_{t-1}-x_{t}}x_{t};
\end{eqnarray}
and
\begin{eqnarray}\label{By3}x_{t-2-i}&=&x_{t-1}+\frac{\lambda(E_{(t-2-i)\backslash (t-1)},\vec{x})}{\lambda(E_{(t-2-i)(t-1)},\vec{x})}\nonumber \\
&=&x_{t-1}+\frac{(x_{t-3-(a-i-2)}+x_{t-3-(a-i-1)}+\cdots+x_{t-3})-x_{t-2-i}}{1-x_{t-2-i}-x_{t-1}-x_{t}}x_{t} \nonumber \\
&\leq& x_{t-1}+\frac{(x_{i+1}+x_{i+2}+\cdots+ x_{a-1})-x_{t-2-i}}{1-x_{t-2-i}-x_{t-1}-x_{t}}x_{t}
\end{eqnarray}
since $a\le t-2$.  Adding (\ref{By2}) and (\ref{By3}), we obtain that
\begin{eqnarray*}
x_{1}+x_{t-2-i}&\leq& x_{t-2}+ x_{t-1}+\frac{(x_{1}+\cdots+ x_{a+1})-x_{t-2-i}}{1-x_{t-2-i}-x_{t-1}-x_{t}}x_{t} \\
&\leq& x_{t-2}+ x_{t-1}+\frac{(x_{1}+\cdots+ x_{t-2})-x_{t-2-i}}{1-x_{t-2-i}-x_{t-1}-x_{t}}x_{t}\\
&=& x_{t-2}+ x_{t-1}+\frac{1-x_{t-1}-x_{t}-x_{t-2-i}}{1-x_{t-2-i}-x_{t-1}-x_{t}}x_{t} \\
&=& x_{t-2}+x_{t-1}+x_{t}.
\end{eqnarray*}
So, (\ref{ineqadd}) is true. This implies that $\lambda(E_{(t-2)(t-1)},\vec{x})\leq\lambda(E_{1(t-2-i)},\vec{x})$. In what follows, we divide the rest of the proof into three cases: $a=2i+1$, $a=2i+2$, and $a \ge 2i+3$.

We first consider the case that $a\geq 2i+3$.
By Remark \ref{r1}(b), we have $x_{1}=x_{2}=\cdots=x_{t-a-2+i}$ and $x_{t-2-i}=\cdots=x_{t-3}$. Hence $\lambda(C_{3,m},\vec{x})-\lambda(G,\vec{x})=i(x_{t-2-i}x_{t-2}x_{t}-x_{1}x_{t-1}x_{t})$. Also by Remark \ref{r1}(b), we have
\begin{eqnarray*}\label{eqII}x_{1}&=&x_{t-2-i}+\frac{\lambda(E_{1\setminus (t-2-i)},\vec{x})}{\lambda(E_{1(t-2-i)},\vec{x})}\nonumber \\
&=&x_{t-2-i}+\frac{(x_{t-1}+x_{t-2})x_{t}}{\lambda(E_{1(t-2-i)},\vec{x})},
\end{eqnarray*}
and
\begin{eqnarray*}x_{t-2}&=&x_{t-1}+\frac{\lambda(E_{(t-2)\setminus (t-1)},\vec{x})}{\lambda(E_{(t-2)(t-1)},\vec{x})}\nonumber \\
&=&x_{t-1}+\frac{(x_{t-3-i}+\cdots +x_{t-a-1+i})x_{t}}{\lambda(E_{(t-2)(t-1)},\vec{x})}.
\end{eqnarray*}

Recall that $a\geq 2i+3$ and $\lambda(E_{(t-2)(t-1)},\vec{x})\leq\lambda(E_{1(t-2-i)},\vec{x})$.
We have $x_{t-2}-x_{t-1}\geq x_{1}-x_{t-2-i}$.
Hence
\begin{eqnarray}\label{eq37}
\lambda(C_{3,m},\vec{x})-\lambda(G,\vec{x})&=&i(x_{t-2-i}x_{t-2}x_{t}-x_{1}x_{t-1}x_{t})\nonumber \\
&=& i[x_{t-2-i}(x_{t-2}+x_{t-1}-x_{t-1})x_{t}-x_{1}x_{t-1}x_{t}] \nonumber \\
&\geq& i[ x_{t-2-i}(x_{t-1}+x_{1}-x_{t-2-i})x_{t}-x_{1}x_{t-1}x_{t}]\nonumber \\
&=&i(x_{t-2-i}-x_{t-1})(x_{1}-x_{t-2-i})x_{t}\nonumber \\
&\geq& 0.
\end{eqnarray}
Therefore $\lambda(C_{3,m})\geq \lambda(C_{3,m},\vec{x})\geq \lambda(G, \vec{x})=\lambda(G)$ in this case.

Next, we consider the case that $a=2i+2$. Let $G'=G\bigcup \{(t-2-i)(t-2)t\}\backslash \{(t-4-i)(t-1)t\}$, then $\lambda(G')\leq \lambda(C_{3,m})$ by the case $a=2(i-1)+4$. (Note that $G'=C_{3,m}$ when $i-1=0$.)  So it is sufficient to prove that $\lambda(G)\leq\lambda(G')$. Clearly,
\begin{equation}\label{e0}
\lambda(G',\vec{x})-\lambda(G,\vec{x})=x_{t-2-i}x_{t-2}x_{t}-x_{t-4-i}x_{t-1}x_{t}=x_{t-2-i}x_{t-2}x_{t}-x_{1}x_{t-1}x_{t}.
\end{equation}
Consider a new weighting ${\vec y}=(y_1, y_2, \ldots, y_t)$ given by $y_j=x_j$ for $j\neq t-4-i$, $j\neq t-2-i$ and $y_{t-4-i}=x_{t-4-i}-\delta$, $y_{ t-2-i}=x_{ t-2-i}+\delta$. Then
\begin{eqnarray*}
\lambda(G',\vec{y})-\lambda(G',\vec{x})&=& \delta[\lambda(E'_{t-2-i},\vec{x})-\lambda(E'_{t-4-i},\vec{x})]-\delta^{2}\lambda(E'_{(t-4-i)(t-2-i)},\vec{x}) \nonumber\\
&=&\delta(x_{t-4-i}-x_{t-2-i})\lambda(E'_{(t-4-i)(t-2-i)},\vec{x})-\delta^{2}\lambda(E'_{(t-4-i)(t-2-i)},\vec{x}).
\end{eqnarray*}
Let $\delta=\frac{x_{t-4-i}-x_{t-2-i}}{2}$. Clearly, ${\vec y}=(y_1, y_2, \ldots, y_t)$ is also a legal weighting for $G$ and
\begin{eqnarray}\label{eq34}
\lambda(G',\vec{y})-\lambda(G',\vec{x})=\frac{(x_{t-4-i}-x_{t-2-i})^{2}}{4}\lambda(E'_{(t-4-i)(t-2-i)},\vec{x})=\frac{(x_{1}-x_{t-2-i})^{2}}{4}\lambda(E_{1(t-2-i)},\vec{x}).
\end{eqnarray}
Let ${\vec z}=(z_1,z_2, \ldots, z_t)$ given by $z_j=y_j$ for $j \neq t-2$, $j \neq t-1$ and $z_{t-2}=y_{t-2}+\eta$, $z_{t-1}=y_{t-1}-\eta$. Then
\begin{eqnarray}\label{eq33}
\lambda(G',\vec{z})-\lambda(G',\vec{y})&=& \eta[\lambda(E'_{t-2},\vec{y})-\lambda(E'_{t-1},\vec{y})]-\eta^{2}\lambda('E_{(t-2)(t-1)},\vec{y}) \nonumber\\
&=&\eta[(y_{t-2-i}y_{t}+y_{t-3-i}y_{t}+y_{t-4-i}y_{t})-(y_{t-2}-y_{t-1})\lambda(E'_{(t-2)(t-1)},\vec{y})]\nonumber\\
& &-\eta^{2}\lambda(E'_{(t-2)(t-1)},\vec{y}).
\end{eqnarray}
Let
\begin{eqnarray*}
\eta &= &\frac{(y_{t-2-i}+y_{t-3-i}+y_{t-4-i})y_{t}-(y_{t-2}-y_{t-1})\lambda(E'_{(t-2)(t-1)},\vec{y})}{2\lambda(E'_{(t-2)(t-1)},\vec{y})} \\
&=& \frac{(x_{t-2-i}+x_{t-3-i}+x_{t-4-i})x_{t}-(x_{t-2}-x_{t-1})
\lambda(E_{(t-2)(t-1)},\vec{x})}{2\lambda(E_{(t-2)(t-1)},\vec{x})}.
\end{eqnarray*}
By Remark \ref{r1}(b), we have
\begin{eqnarray}\label{eq35}
x_{t-2}=x_{t-1}+\frac{x_{t-3-i}x_{t}}{\lambda(E_{(t-2)(t-1)},\vec{x})}.
\end{eqnarray}
Hence, $\eta=\frac{(x_{t-2-i}+x_{t-4-i})x_{t}}{2\lambda(E_{(t-2)(t-1)},\vec{x})}$ and ${\vec z}=(z_1, z_2, \ldots, z_t)$ is also a legal weighting for $G$ and
\begin{eqnarray}\label{eq344}
\lambda(G',\vec{z})-\lambda(G',\vec{y})=\frac{(x_{t-2-i}+x_{t-4-i})^{2}x_{t}^{2}}{4\lambda(E_{(t-2)(t-1)},\vec{x})}.
\end{eqnarray}
By Remark \ref{r1}(b), we have
\begin{eqnarray}\label{eq35}
x_{1}=x_{t-i-2}+\frac{x_{t-2}x_{t}+x_{t-1}x_{t}}{\lambda(E_{1(t-i-2)},\vec{x})}.
\end{eqnarray}
Combing (\ref{e0}), (\ref{eq34}), (\ref{eq344}), and (\ref{eq35}), we have
\begin{eqnarray*}
\lambda(G',\vec{z})-\lambda(G,\vec{x})&=&x_{t-2-i}x_{t-2}x_{t}-x_{1}x_{t-1}x_{t}+\frac{(x_{t-2}+x_{t-1})^{2}x_{t}^{2}}{4\lambda(E_{1(t-i-2)},\vec{x})}+\frac{(x_{t-2-i}+x_{t-4-i})^{2}x_{t}^{2}}{4\lambda(E_{(t-2)(t-1)},\vec{x})}\\
&=&[x_{1}-\frac{x_{t-2}x_{t}+x_{t-1}x_{t}}{\lambda(E_{1(t-i-2)},\vec{x})}]x_{t-2}x_{t}
-x_{1}x_{t-1}x_{t} \\
&&+\frac{(x_{t-2}+x_{t-1})^{2}x_{t}^{2}}{4\lambda(E_{1(t-i-2)},\vec{x})}+\frac{(x_{t-2-i}+x_{t-4-i})^{2}x_{t}^{2}}{4\lambda(E_{(t-2)(t-1)},\vec{x})} \\&\ge & {x_t^2 \over 4\lambda(E_{1(t-i-2)},\vec{x})}[-4(x_{t-2}+x_{t-1})x_{t-2}+(x_{t-2}+x_{t-1})^{2}+4x_{t-2}^2)]\\
&=& {x_t^2 \over 4\lambda(E_{1(t-i-2)},\vec{x})}(x_{t-2}-x_{t-1})^{2}\\&\ge& 0.
\end{eqnarray*}
Hence $\lambda(G')\geq\lambda(G',\vec{z})\geq\lambda(G,\vec{x})=\lambda(G)$ in this case.

What remains is the case that $a=2i+1$. Let $G''=G\bigcup \{(t-2-i)(t-2)t\}\backslash \{(t-3-i)(t-1)t\}$, then $\lambda(G'')\leq \lambda(C_{3,m})$ by the case $a=2(i-1)+3$. (Note that $G'=C_{3,m}$ when $i-1=0$.)  So it is sufficient to prove that $\lambda(G)\leq\lambda(G'')$. Clearly,
\begin{equation}\label{e00}
\lambda(G'',\vec{x})-\lambda(G,\vec{x})=x_{t-2-i}x_{t-2}x_{t}-x_{t-3-i}x_{t-1}x_{t}=x_{t-2-i}x_{t-2}x_{t}-x_{1}x_{t-1}x_{t}.
\end{equation}
Consider a new weighting ${\vec u}=(u_1, u_2, \ldots, u_t)$ given by $u_j=x_j$ for $j \neq t-2-i$, $j \neq t-3-i$ and $u_{t-2-i}=x_{t-2-i}+\alpha$, $u_{t-3-i}=x_{t-3-i}-\alpha$. Then
\begin{eqnarray*}
\lambda(G'',\vec{u})-\lambda(G'',\vec{x})&=& \alpha(x_{t-3-i}-x_{t-2-i})\lambda(E''_{(t-3-i)(t-2-i)},\vec{x})-\alpha^{2}\lambda(E''_{(t-3-i)(t-2-i)},\vec{x}).
\end{eqnarray*}
Let $\alpha=\frac{x_{t-3-i}-x_{t-2-i}}{2}$. Clearly, ${\vec u}=(u_1, u_2, \ldots, u_t)$ is also a legal weighting and
\begin{eqnarray}\label{eq51}
\lambda(G'',\vec{u})-\lambda(G'',\vec{x})=\frac{(x_{t-i-3}-x_{t-2-i})^{2}}{4}\lambda(E''_{(t-3-i)(t-2-i)},\vec{x})=\frac{(x_{1}-x_{t-2-i})^{2}}{4}\lambda(E_{1(t-2-i)},\vec{x}).
\end{eqnarray}
Let ${\vec v}=(v_1,v_2, \ldots, v_t)$ given by $v_j=u_j$ for $j \neq t-2$, $j \neq t-1$ and $v_{t-2}=u_{t-2}+\beta$, $v_{t-1}=u_{t-1}-\beta$. Then
\begin{eqnarray}\label{eq33}
\lambda(G'',\vec{v})-\lambda(G'',\vec{u})&=& \beta[\lambda(E''_{t-2},\vec{u})-\lambda(E''_{t-1},\vec{u})]-\beta^{2}\lambda(E''_{(t-2)(t-1)},\vec{u}) \nonumber\\
&=&\beta(u_{t-2-i}u_{t}+u_{t-3-i}u_{t})-\beta^{2}\lambda(E''_{(t-2)(t-1)},\vec{u}).
\end{eqnarray}
Let $\beta=\frac{u_{t-2-i}u_{t}+u_{t-3-i}u_{t}}{2\lambda(E''_{(t-2)(t-1)},\vec{u})}$. Clearly, $\beta<u_{t}$. Hence, ${\vec v}=(v_1, v_2, \ldots, v_t)$ is also a legal weighting for $G$ and
\begin{eqnarray}\label{eq52}
\lambda(G'',\vec{v})-\lambda(G'',\vec{u})=\frac{(u_{t-2-i}+u_{t-3-i})^{2}u_{t}^{2}}{4\lambda(E''_{(t-2)(t-1)},\vec{u})}=\frac{(x_{t-2-i}+x_{t-3-i})^{2}x_{t}^{2}}{4\lambda(E_{(t-2)(t-1)},\vec{x})}.
\end{eqnarray}
By Remark \ref{r1}(b), we have $x_{t-2}=x_{t-1}$ and
\begin{eqnarray}\label{eq53}
x_{1}=x_{t-2-i}+\frac{2x_{t-1}x_{t}}{\lambda(E_{1(t-2-i)},\vec{x})}.
\end{eqnarray}
Combing (\ref{e00}), (\ref{eq51}), (\ref{eq52}), and (\ref{eq53}), we have
\begin{eqnarray*}
\lambda(G'',\vec{v})-\lambda(G,\vec{x})&=&x_{t-2-i}x_{t-2}x_{t}-x_{1}x_{t-1}x_{t}+\frac{(x_{t-2-i}+x_{t-3-i})^{2}x_{t}^{2}}{4\lambda(E_{(t-2)(t-1)},\vec{x})}
+\frac{x_{t-1}^{2}x_{t}^{2}}{\lambda(E_{1(t-2-i)},\vec{x})} \\
&=&-\frac{2x_{t-1}^{2}x_{t}^{2}}{\lambda(E_{1(t-2-i)},\vec{x})}
+\frac{(x_{t-2-i}+x_{t-3-i})^{2}x_{t}^{2}}{4\lambda(E_{(t-2)(t-1)},\vec{x})}
+\frac{x_{t-1}^{2}x_{t}^{2}}{\lambda(E_{1(t-2-i)},\vec{x})} \\
&\ge& 0
\end{eqnarray*}
since $\lambda(E_{(t-2)(t-1)},\vec{x})\leq \lambda(E_{1(t-2-i)},\vec{x})$.   Hence $\lambda(G'')\geq\lambda(G'',\vec{z})\geq\lambda(G,\vec{x})=\lambda(G).$
This completes the proof of Theorem \ref{Lemma0}. \qed

\subsection{Proof of Theorem \ref{mainresult}}

\begin{remark}\label{reducetot2}
Let $G$ be a left-compressed extremal
$3$-graph with $m$ edges satisfying $|E(G) \Delta E(C_{3,m})| \le 6$.  Let $t$ be a positive integer such that ${t-1 \choose 3} \le m < {t \choose 3}$. To show $\lambda(G)\leq\lambda(C_{3,m})$, we can assume $G$ is  on  vertex set $[t]$.
\end{remark}

\noindent{\em Proof.} The proof is excatly the same as the proof of Remark \ref{reducetot}.\qed

 The partial ordered diagram (Figure 1) on all triples on  $[t]$ as described below is useful to help us to analyze all possible left-compressed 3-graphs on $[t]$ systematically

An  $r$-tuple  $i_1 i_2 \cdots i_r$ is called a {\it descendant  } of an  $r$-tuple  $j_1j_2 \cdots j_r$ if $i_s\le j_s$ for each $1\le s\le r$, and $i_1+i_2+\cdots +i_r < j_1+j_2+\cdots +j_r$. In this case, $j_1j_2\cdots j_r$  is called an {\it ancestor } of $i_1 i_2\cdots i_r$.  The  $r$-tuple $i_1 i_2 \cdots i_r$   is called a {\it direct descendant} of $j_1j_2 \cdots j_r$ if  $i_1 i_2 \cdots i_r$ is a descendant of $j_1j_2 \cdots j_r$ and $j_1+j_2+\cdots +j_r=i_1+i_2+\cdots +i_r +1$.  We say that $j_1 j_2\cdots j_r$ has lower hierarchy than $i_1i_2\cdots i_r$ if $j_1 j_2 \cdots j_r$  is  an ancestor of $i_1i_2\cdots i_r$. This is a partial order on the set of all $r$-tuples.  Figure 1 is a Hessian diagram on all triples on vertex set $[t]$. In this diagram, $i_1 i_2 i_3$ and  $j_1j_2 j_3$ are connected by an edge if and only if $i_1 i_2 i_3$   is  a  direct descendant of $j_1j_2 j_3$.

\begin{figure}[!htb]
\centering
\includegraphics{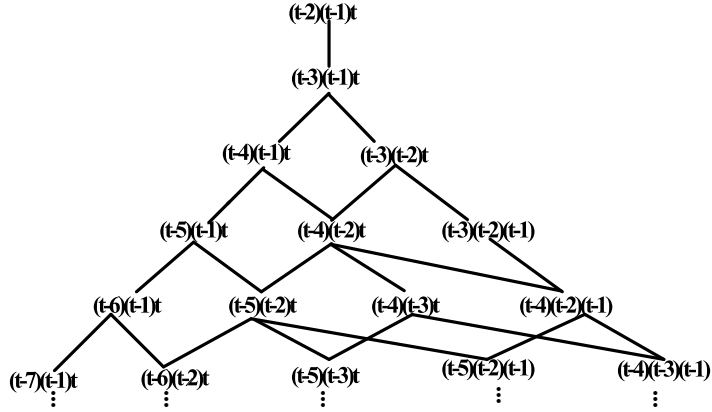}
\caption{Hessian Diagram on $[t]^{(3)}$}
\end{figure}

Let us be aware of the following simple observation for left-compressed $r$-graphs.
\begin{remark}\label{releftcom}
An $r$-graph $G$ is left-compressed if and only if all descendants of an edge of $G$ are edges of $G$. Equivalently, if an $r$-tuple is not an edge of $G$, then none of its ancestors will be an edge of $G$.
\end{remark}

By analyzing possible cases under the assumption that $G$ is a left-compressed extremal graph on $[t]$ satisfying $|E(G) \triangle E(C_{3,m})|\leq 6$, we give several lemmas to cover
 the possible cases below.
% The first one was proved in  \cite{TPZZ3}.

%\begin{lemma}(Tang,Peng,Zhang and Zhao \cite{TPZZ3})\label{Lemmatpzz2} Let $m$ and $t$ be positive integers satisfying $${t-1 \choose 3} \le m \le {t-1 \choose 3} + {t-2 \choose 2}.$$  Let $G=([t],E)$ be a left-compressed 3-graph  with $m$ edges satisfying $|E_{(t-1)t}|\le3$. Then $\lambda (G)\leq \lambda ([t-1]^{(3)})$.
%\end{lemma}

%\begin{center}
%\includegraphics{figure1_tang}
%\end{center}
%\begin{figure}[!htb]
%\centering
%\includegraphics{figure1_tang}
%\caption{}
%\end{figure}

%In the following lemmas, we assume that $a$ is an integer satisfying $3\leq a\leq t-3$.

Using Theorem \ref{Lemma0}, we deal with the case when the graph $G$ contains a clique of order $t-1$ in Lemma \ref{Lemma1}.

\begin{lemma} \label{Lemma1} Let $m$, $a$ and $t$ be  positive integers satisfying  $m={t \choose 3}-a$, where $3\le a\le t-2$. Let $G=([t],E)$ be a left-compressed $3$-graph  with $m$ edges and a clique of order  $t-1$. If $|E\triangle E(C_{3,m})|\leq 6$, then $\lambda(G)\leq\lambda(C_{3,m})$.
\end{lemma}

\noindent{\em Proof.} If the triple with the  minimum colex ordering in $G^{c}$ is $(t-2-i)(t-2)t$, where $i=1,2,3$. Then $\lambda(G)\leq \lambda (C_{3,m})$ by Theorem \ref{Lemma0}. So  we can assume that the triple with the minimum colex ordering in $G^{c}$ is $(t-4)(t-3)t$. Let $\vec{x}=(x_{1},x_{2},\ldots ,x_{t})$  be an optimal weighting for $G$ satisfying $x_1 \ge x_2 \ge \ldots \ge x_t \ge 0$.

First we consider the case $a\geq 7$. Let $G'=G\bigcup\{(t-4)(t-3)t\}\backslash \{(t-5)(t-2)t\}$, then $\lambda(G')\leq \lambda(C_{3,m})$ by Theorem \ref{Lemma0}. So it is sufficient to prove that $\lambda(G)\leq \lambda(G')$. Note that
\begin{eqnarray}\label{eq80}
\lambda(G',\vec{x})-\lambda(G,\vec{x})=x_{t-4}x_{t-3}x_{t}-x_{t-5}x_{t-2}x_{t}.
\end{eqnarray}
Consider a new weighting ${\vec y}=(y_1, y_2, \ldots, y_t)$ given by $y_j=x_j$ for $j \neq t-5$, $j \neq t-4$ and $y_{t-5}=x_{t-5}-\delta$, $y_{t-4}=x_{t-4}+\delta$.
Therefore
\begin{eqnarray*}
\lambda(G',\vec{y})-\lambda(G',\vec{x})&=& \delta[\lambda(E'_{t-4},\vec{x})-\lambda(E'_{t-5},\vec{x})]-\delta^{2}\lambda(E'_{(t-5)(t-4)},\vec{x}) \nonumber\\
&=&\delta(x_{t-5}-x_{t-4})\lambda(E'_{(t-5)(t-4)},\vec{x})-\delta^{2}\lambda(E'_{(t-5)(t-4)},\vec{x}).
\end{eqnarray*}
Let $\delta=\frac{x_{t-5}-x_{t-4}}{2}$. Clearly, ${\vec y}=(y_1, y_2, \ldots, y_t)$ is also a legal weighting for $G$.  Note that $\lambda(E'_{(t-5)(t-4)},\vec{x})=\lambda(E_{(t-5)(t-4)},\vec{x})$, so
\begin{eqnarray}\label{eq81}
\lambda(G',\vec{y})-\lambda(G',\vec{x})=\frac{(x_{t-5}-x_{t-4})^{2}}{4}\lambda(E'_{(t-5)(t-4)},\vec{x})=\frac{(x_{t-5}-x_{t-4})^{2}}{4}\lambda(E_{(t-5)(t-4)},\vec{x}).
\end{eqnarray}
Let ${\vec z}=(z_1,z_2, \ldots, z_t)$ given by $z_j=y_j$ for $j \neq t-3$, $j \neq t-2$ and $z_{t-3}=y_{t-3}+\eta$, $z_{t-2}=y_{t-2}-\eta$. Then
\begin{eqnarray}
\lambda(G',\vec{z})-\lambda(G',\vec{y})&=& \eta[\lambda(E'_{t-3},\vec{y})-\lambda(E'_{t-2},\vec{y})]-\eta^{2}\lambda(E'_{(t-3)(t-2)},\vec{y}) \nonumber\\
&=&\eta(y_{t-4}y_{t}+y_{t-5}y_{t})-\eta^{2}\lambda(E'_{(t-3)(t-2)},\vec{y}).
\end{eqnarray}
Let $\eta=\frac{y_{t-5}y_{t}+y_{t-4}y_{t}}{2\lambda(E'_{(t-3)(t-2)},\vec{y})}$. Clearly, $\eta<y_{t}$. Hence, ${\vec z}=(z_1, z_2, \ldots, z_t)$ is also a legal weighting and
\begin{eqnarray}\label{eq82}
\lambda(G',\vec{z})-\lambda(G',\vec{y})=\frac{(y_{t-5}+y_{t-4})^{2}y_{t}^{2}}{4\lambda(E'_{(t-3)(t-2)},\vec{y})}=\frac{(x_{t-5}+x_{t-4})^{2}x_{t}^{2}}{4\lambda(E_{(t-3)(t-2)},\vec{x})}.
\end{eqnarray}
By Remark \ref{r1}(b), we have $x_{t-4}=x_{t-3}=x_{t-2}$ and
\begin{eqnarray}\label{eq83}
x_{t-5}=x_{t-4}+\frac{2x_{t-2}x_{t}}{\lambda(E_{(t-5)(t-4)},\vec{x})}.
\end{eqnarray}
 Also,
\begin{eqnarray}\label{eq84}
\lambda(E_{(t-5)(t-4)},\vec{x})-\lambda(E_{(t-3)(t-2)},\vec{x})\geq x_{t-2}+x_{t}-x_{t-5}\geq 0
\end{eqnarray}
since $x_{t-5}\leq x_{1}=x_{t-1}+\frac{\lambda(E_{1\backslash (t-1)},\vec{x})}{\lambda(E_{1(t-1)},\vec{x})}\leq x_{t-1}+\frac{(x_{2}+\cdots+x_{t-2})x_{t}}{x_{2}+\cdots+x_{t-2}+x_{t}}\leq x_{t-1}+x_{t}\leq x_{t-2}+x_{t} $.

Combining (\ref{eq80}), (\ref{eq81}), (\ref{eq82}) and (\ref{eq83}), we have
\begin{eqnarray*}
\lambda(G',\vec{z})-\lambda(G,\vec{x})&\ge &x_{t-4}x_{t-3}x_{t}-x_{t-5}x_{t-2}x_{t}+\frac{(x_{t-5}+x_{t-4})^{2}x^{2}_{t}}{4\lambda(E_{(t-3)(t-2)},\vec{x})}
+\frac{x_{t-2}^{2}x_{t}^{2}}{\lambda(E_{(t-5)(t-4)},\vec{x})}  \nonumber\\
&=&-\frac{2x_{t-2}^{2}x_{t}^{2}}{\lambda(E_{(t-5)(t-4)},\vec{x})}
+\frac{(x_{t-5}+x_{t-4})^{2}x^{2}_{t}}{4\lambda(E_{(t-3)(t-2)},\vec{x})}
+\frac{x_{t-2}^{2}x_{t}^{2}}{\lambda(E_{(t-5)(t-4)},\vec{x})}\\&\geq& 0
\end{eqnarray*}
since $\lambda(E_{(t-5)(t-4)},\vec{x})\ge \lambda(E_{(t-3)(t-2)},\vec{x})$ in view of (\ref{eq84}). Hence $\lambda(G')\geq\lambda(G',\vec{z})\geq\lambda(G,\vec{x})=\lambda(G)$.

Next, we consider the case when $a \le 6$. Clearly the lemma holds if $a=3$, $4$ and $5$ in view of Theorem \ref{Lemma0}.
So the only remaining case is that $a=6$. In view of Figure 1, we have $E(G)=E=[t]^{(3)}\setminus\{(t-2)(t-1)t,(t-3)(t-1)t,(t-4)(t-1)t,(t-3)(t-2)t,(t-4)(t-2)t,(t-4)(t-3)t\}$. In this case, $E(C_{3,m})=E''=[t]^{(3)}\setminus\{(t-2)(t-1)t,(t-3)(t-1)t,(t-4)(t-1)t,(t-5)(t-1)t,(t-6)(t-1)t,(t-7)(t-1)t\}$. By Remark \ref{r1}(b), we have
$$x_{1}=x_{2}=\cdots=x_{t-5}=\alpha;\quad x_{t-4}=x_{t-3}=x_{t-2}=x_{t-1}=\beta; \quad x_{t}=\gamma; $$
and
\begin{equation}\label{eq000}
\alpha-\beta=x_{t-5}-x_{t-4}=\frac{\lambda(E_{(t-5)\backslash (t-4)},\vec{x})}{\lambda(E_{(t-5)(t-4)},\vec{x})}={ x_{t-2}x_{t}+x_{t-3}x_{t}+x_{t-4}x_{t} \over \lambda(E_{(t-5)(t-4)},\vec{x})}={3\beta\gamma \over \lambda(E_{(t-5)(t-4)},\vec{x})}\le \gamma.
\end{equation}
Then
\begin{equation}\label{eq0000}
\lambda(C_{3,m},\vec{x})-\lambda(G,\vec{x})=3\beta^2\gamma-3\alpha\beta\gamma=-\frac{9\beta^2\gamma^2}{\lambda(E_{(t-5)(t-4)},\vec{x})}=
-\frac{9\beta^2\gamma^2}{\lambda(E_{(t-6)(t-4)},\vec{x})}=-\frac{9\beta^2\gamma^2}{\lambda(E_{(t-7)(t-4)},\vec{x})}.
\end{equation}
 Let ${\vec u}=(u_1,u_2, \ldots, u_t)$ given by $u_j=x_j$ for $j\neq t-7$, $j\neq t-4$ and $u_{t-7}=x_{t-7}-\delta$, $u_{t-4}=x_{t-4}+\delta$. Then
\begin{eqnarray*}
\lambda(C_{3,m},\vec{u})-\lambda(C_{3,m},\vec{x})&=& \delta[\lambda(E''_{t-4},\vec{x})-\lambda(E''_{t-7},\vec{x})]-\delta^{2}\lambda(E''_{(t-7)(t-4)},\vec{x}) \nonumber\\
&=&\delta(x_{t-7}-x_{t-4})\lambda(E''_{(t-7)(t-4)},\vec{x})-\delta^{2}\lambda(E''_{(t-7)(t-4)},\vec{x}).
\end{eqnarray*}
Let $\delta=\frac{x_{t-7}-x_{t-4}}{2}=\frac{\alpha-\beta}{2}$. Clearly, ${\vec u}=(u_1, u_2, \ldots, u_t)$ is also a legal weighting for $G$, $\lambda(E''_{(t-7)(t-4)},\vec{x})=\lambda(E_{(t-7)(t-4)},\vec{x})$ and
\begin{eqnarray}\label{eq861}
\lambda(C_{3,m},\vec{u})-\lambda(C_{3,m},\vec{x})&=&\frac{(\alpha-\beta)^{2}}{4}\lambda(E''_{(t-7)(t-4)},\vec{x})\nonumber\\
&=&\frac{9\beta^2\gamma^2}{4[\lambda(E_{(t-7)(t-4)},\vec{x})]^{2}}\lambda(E''_{(t-7)(t-4)},\vec{x})\nonumber\\&=&\frac{9\beta^2\gamma^2}{4\lambda(E_{(t-7)(t-4)},\vec{x})}.
\end{eqnarray}
Similarly, let ${\vec v}=(v_1,v_2, \ldots, v_t)$ given by $v_j=u_j$ for $j\neq t-6$, $j\neq t-3$ and $v_{t-6}=u_{t-6}-\frac{\alpha-\beta}{2}$, $v_{t-3}=u_{t-3}+\frac{\alpha-\beta}{2}$. Then
\begin{eqnarray}\label{eq862}
\lambda(C_{3,m},\vec{v})-\lambda(C_{3,m},\vec{u})=\frac{9\beta^2\gamma^2}{4\lambda(E_{(t-7)(t-4)},\vec{x})}.
\end{eqnarray}
Let ${\vec w}=(w_1,w_2, \ldots, w_t)$ given by $w_j=v_j$ for $j\neq t-5$, $j\neq t-2$ and $w_{t-5}=v_{t-5}-\frac{\alpha-\beta}{2}$, $w_{t-2}=v_{t-2}+\frac{\alpha-\beta}{2}$. Then
\begin{eqnarray}\label{eq3863}
\lambda(C_{3,m},\vec{w})-\lambda(C_{3,m},\vec{v})=\frac{9\beta^2\gamma^2}{4\lambda(E_{(t-7)(t-4)},\vec{x})}.
\end{eqnarray}
Let ${\vec p}=(p_1,p_2, \ldots, p_t)$ given by $p_j=w_j$ for $j\neq t-1$, $j\neq t$ and $p_{t-1}=w_{t-1}-\frac{\beta-\gamma}{2}$, $p_{t}=w_{t}+\frac{\beta-\gamma}{2}$. Then
\begin{eqnarray}\label{eq864}
\lambda(C_{3,m},\vec{p})-\lambda(C_{3,m},\vec{w})=\frac{9\beta^4}{4\lambda(E_{(t-1)t},\vec{x})}.
\end{eqnarray}
Note that
\begin{eqnarray}\label{eq865}
 \lambda(E_{(t-5)(t-4)},\vec{x})-\lambda(E_{(t-1)t},\vec{x})=3\beta+\gamma-\alpha\geq 0
\end{eqnarray}
by (\ref{eq000}).

Combining (\ref{eq0000}),  (\ref{eq861}), (\ref{eq862}), (\ref{eq3863}), (\ref{eq864}), and (\ref{eq865}), we have
\begin{eqnarray}
\lambda(C_{3,m},\vec{p})-\lambda(G,\vec{x})=\frac{9\beta^4}{4\lambda(E_{(t-1)t},\vec{x})}-\frac{9\beta^2\gamma^2}{4\lambda(E_{(t-5)(t-4)},\vec{x})}\geq 0.
\end{eqnarray}
Hence $\lambda(C_{3,m})\geq\lambda(C_{3,m},\vec{p})\geq\lambda(G,\vec{x})=\lambda(G)$. This completes the proof of Lemma \ref{Lemma1}.\qed

Using Lemma \ref{Lemma1}, we prove the next four lemmas which cover the cases when the 3-graph $G$ does not contain a clique of order $t-1$.

\begin{lemma} \label{Lemma2} Let $G$ and $G'$ be  left-compressed $3$-graphs on vertex set $[t]$ with $m={t \choose 3}-a$ edges, where $5\le a\le t-2$, satisfying $|E(G) \triangle E(C_{3,m})|=|E(G') \triangle E(C_{3,m})|=4$ and the triples with the minimum colex ordering in $G^{c}$ and $G'^{c}$ are $(t-3)(t-2)(t-1)$ and $(t-4)(t-2)t$ respectively. Then $ \lambda(G)\leq \lambda(G')\leq \lambda (C_{3,m})$.
\end{lemma}
{\em Proof.} By Lemma \ref{Lemma1}, $\lambda(G')\leq \lambda (C_{3,m})$. So it is sufficient to show $ \lambda(G)\leq \lambda(G')$. Let $\vec{x}=(x_{1},x_{2},\ldots ,x_{t})$  be an optimal weighting for $G$ satisfying $x_1 \ge x_2  \geq\cdots\geq x_{t}\geq 0$.
By Remark \ref{r1}(b), $x_{t-2}=x_{t-3}$ and $ x_{t-1}=x_{t}$. Hence
\begin{eqnarray}\label{eq91}
\lambda(G',\vec{x})-\lambda(G,\vec{x})=(x_{t-3}-x_{t-4})x_{t-2}x_{t-1}.
\end{eqnarray}
Consider a new weighting ${\vec y}=(y_1, y_2, \ldots, y_t)$ given by $y_j=x_j$ for $j\neq t-4$, $j\neq t-3$ and $y_{t-4}=x_{t-4}-\delta$, $y_{t-3}=x_{t-3}+\delta$. Then
\begin{eqnarray}
\lambda(G',\vec{y})-\lambda(G',\vec{x})&=& \delta[\lambda(E'_{t-3},\vec{x})-\lambda(E'_{t-4},\vec{x})]-\delta^{2}\lambda(E'_{(t-4)(t-3)},\vec{y}) \nonumber\\
&=&\delta(x_{t-4}-x_{t-3})\lambda(E'_{(t-4)(t-3)},\vec{x})-\delta^{2}\lambda(E'_{(t-4)(t-3)},\vec{x}).
\end{eqnarray}
Let $\delta=\frac{x_{t-4}-x_{t-3}}{2}$. Clearly, ${\vec y}=(y_1, y_2, \ldots, y_t)$ is also a legal weighting for $G$. Also note that $\lambda(E'_{(t-4)(t-3)},\vec{x})=\lambda(E_{(t-4)(t-3)},\vec{x})$. Hence
\begin{eqnarray}\label{eq92}
\lambda(G',\vec{y})-\lambda(G',\vec{x})=\frac{(x_{t-4}-x_{t-3})^{2}}{4}\lambda(E'_{(t-4)(t-3)},\vec{x})=\frac{(x_{t-4}-x_{t-3})^{2}}{4}\lambda(E_{(t-4)(t-3)},\vec{x}).
\end{eqnarray}
Let ${\vec z}=(z_1, z_2, \ldots, z_t)$ be given by $z_i=y_i$ for $i\neq t-1$, $i\neq t$ and $z_{t-1}=y_{t-1}+\eta$, $z_{t}=y_{t}-\eta$. Then
\begin{eqnarray}
\lambda(G',\vec{z})-\lambda(G',\vec{y})&=& \eta[\lambda(E'_{t-1},\vec{y})-\lambda(E'_{t},\vec{y})]-\eta^{2}\lambda(E'_{(t-1)t},\vec{y}) \nonumber\\
&=&\eta(y_{t-3}y_{t-2}+y_{t-4}y_{t-2})-\eta^{2}\lambda(E'_{(t-1)t},\vec{y})\nonumber\\
&=&\eta(x_{t-3}x_{t-2}+x_{t-4}x_{t-2})-\eta^{2}\lambda(E_{(t-1)t},\vec{x})
\end{eqnarray}
in view of $y_{t-4}+y_{t-3}=x_{t-4}+x_{t-3}$, $y_{t-2}=x_{t-2}$ and $\lambda(E'_{(t-1)t},\vec{y})=\lambda(E_{(t-1)t},\vec{x})$.
Let $\eta=\frac{x_{t-3}x_{t-2}+x_{t-4}x_{t-2}}{2\lambda(E_{(t-1)t},\vec{x})}$. By the condition of $|E(G) \triangle E(C_{3,m})|=4$ we have $\{1,2\} \subseteq E_{(t-1)t}$, so $$\eta\le \frac{x_{t-2}}{2}.$$ Applying Remark \ref{r1}(b), we have
\begin{eqnarray}
x_{t-2}&=&x_{t}+{\lambda(E_{(t-2)\backslash t},\vec{x})\over \lambda(E_{(t-2)t},\vec{x})}\nonumber\\
&\le &x_{t}+\frac{(x_{t-4}+\cdots +x_3)x_{t-1}}{1-x_t-x_{t-1}-x_{t-2}-x_{t-3}}\nonumber\\
&\leq&x_{t}+x_{t-1}\nonumber\\
&=&2x_{t}.
\end{eqnarray}
So $\frac{x_{t-2}}{2}\leq x_{t}$. Recall that  $\eta\leq\frac{x_{t-2}}{2}$. Therefore, $\eta\leq x_{t}$.
Hence, ${\vec z}=(z_1, z_2, \ldots, z_t)$ is also a legal weighting for $G'$, and
\begin{eqnarray}\label{eq93}
\lambda(G',\vec{z})-\lambda(G',\vec{y})=\frac{(x_{t-4}+x_{t-3})^{2}x^{2}_{t-2}}{4\lambda(E_{(t-1)t},\vec{x})}.
\end{eqnarray}
By Remark \ref{r1}(b), we have
\begin{eqnarray}\label{eq94}
x_{t-4}=x_{t-3}+\frac{2x_{t-2}x_{t}}{\lambda(E_{(t-4)(t-3)},\vec{x})}.
\end{eqnarray}
In addition,
\begin{eqnarray}\label{eq95}
\lambda(E_{(t-4)(t-3)},\vec{x})-\lambda(E_{(t-1)t},\vec{x})\geq (1-x_{t-4}-x_{t-3})-(1-x_{t-4}-x_{t-3}-x_{t-2}-x_{t-1}-x_{t})> 0.
\end{eqnarray}
Combing (\ref{eq91}), (\ref{eq92}), (\ref{eq93}), (\ref{eq94}), and (\ref{eq95}), we have
\begin{eqnarray}
\lambda(G',\vec{z})-\lambda(G,\vec{x})&=&\frac{(x_{t-4}-x_{t-3})^{2}\lambda(E_{(t-4)(t-3)},\vec{x})}{4} + \frac{(x_{t-4}+x_{t-3})^{2}x^{2}_{t-2}}{4\lambda(E_{(t-1)t},\vec{x})} \nonumber\\
&& -\frac{2x_{t-2}^{2}x_{t}^{2}}{\lambda(E_{(t-4)(t-3)},\vec{x})}\nonumber\\
&=&\frac{x_{t-2}^{2}x_{t}^2}{\lambda(E_{(t-4)(t-3)},\vec{x})}+\frac{(x_{t-4}+x_{t-3})^{2}x^{2}_{t-2}}{4\lambda(E_{(t-1)t},\vec{x})} -\frac{2x_{t-2}^{2}x_{t}^{2}}{\lambda(E_{(t-4)(t-3)},\vec{x})}\nonumber\\
&\ge& 0.
\end{eqnarray}
Hence $\lambda(G)=\lambda(G,\vec{x})\leq \lambda(G',\vec{z})\leq\lambda(G').$  \qed

\begin{lemma} \label{Lemma2+} Let $G$  be  the left-compressed $3$-graphs on vertex set $[t]$ with $m={t \choose 3}-4$ edges and
$G^{c}=\{(t-2)(t-1)t,(t-3)(t-1)t, (t-3)(t-2)t, (t-3)(t-2)(t-1)\}$.
Then $ \lambda(G)\leq \lambda (C_{3,m})$.
\end{lemma}
{\em Proof.} Let $\vec{x}=(x_{1},x_{2},\ldots ,x_{t})$  be an
optimal weighting for $G$ satisfying $x_1 \ge x_2 \ge \ldots \ge x_t
\ge 0$. By Remark \ref{r1}(b), we have
$$x_1=x_2=\cdots=x_{t-4}=\alpha, \quad x_{t-3}=x_{t-2}=x_{t-1}=x_{t}=\beta,$$
and $$\alpha-\beta=x_{t-4}-x_{t-3}={\lambda(E_{(t-4)\setminus
(t-3)},\vec{x}) \over
\lambda(E_{(t-4)(t-3)},\vec{x})}=\frac{3\beta^2}{\lambda(E_{(t-4)(t-3)},\vec{x})}\le
\beta.$$
 Therefore,
\begin{equation}\label{ale2b}
\alpha\le 2\beta.
\end{equation}
Note that
$E(C_{3,m})=E'=[t]^{(3)}\backslash \{(t-2)(t-1)t,(t-3)(t-1)t, (t-4)(t-1)t,
(t-5)(t-1)t\}$, so
 \begin{eqnarray}\label{g13x}
\lambda(C_{3,m},\vec{x})-\lambda(G,\vec{x})=(x_{t-3}x_{t-2}x_{t-1}+x_{t-3}x_{t-2}x_{t})-(x_{t-4}x_{t-1}x_{t}+x_{t-5}x_{t-1}x_{t})=2(\beta-\alpha)\beta^2.
\end{eqnarray}
Let $\vec {y}=(y_1, y_2, \ldots,
y_t)$ be given by $y_i=x_i$ for $i\neq t-3$, $i\neq t$ and
$y_{t-3}=x_{t-3}+\frac{\alpha-\beta}{3}=\frac{\alpha+2\beta}{3}$,
$y_{t}=x_{t}-\frac{\alpha-\beta}{3}=\frac{4\beta-\alpha}{3}$.
Clearly, $\vec{ y}=(y_1, y_2, \ldots, y_t)$ is a legal weighting for
$C_{3,m}$, and
\begin{eqnarray}\label{eq18e}
\lambda(C_{3,m},\vec{y})-\lambda(C_{3,m},\vec{x})&=&\frac{\alpha-\beta}{3}[\lambda(E'_{t-3},\vec{x})-\lambda(E'_{t},\vec{x})]-
(\frac{\alpha-\beta}{3})^{2}\lambda(E'_{(t-3)t},\vec{x})\nonumber\\
&=&\frac{\alpha-\beta}{3}(x_{t-2}x_{t-1}+x_{t-4}x_{t-1}+x_{t-5}x_{t-1})- (\frac{\alpha-\beta}{3})^{2}\lambda(E'_{(t-3)t},\vec{x})\nonumber\\
&=&\frac{\alpha-\beta}{3}(\beta^2+2\alpha\beta)-\frac{(\alpha-\beta)^2}{9}\lambda(E'_{(t-3)t},\vec{x}).
\end{eqnarray}
Let $\vec{z}=(z_1, z_2, \ldots,
z_t)$ be given by $z_i=y_i$ for $i\neq t-2, i\neq t-1$ and
$z_{t-2}=y_{t-2}+\frac{\alpha-\beta}{3}=\frac{\alpha+2\beta}{3}$,
$z_{t-1}=y_{t-1}-\frac{\alpha-\beta}{3}=\frac{4\beta-\alpha}{3}$.
Clearly $\vec{ z}=(z_1, z_2, \ldots, z_t)$ is also a legal weighting
for $C_{3,m}$, and
\begin{eqnarray}\label{eq19e}
\lambda(C_{3,m},\vec{z})-\lambda(C_{3,m},\vec{y})&=&
\frac{\alpha-\beta}{3}[\lambda(E'_{t-2},\vec{y})-\lambda(E'_{t-1},\vec{y})]-
(\frac{\alpha-\beta}{3})^{2}\lambda(E'_{(t-2)(t-1)},\vec{y})\nonumber\\
&=&\frac{\alpha-\beta}{3}(y_{t-3}y_{t}+y_{t-4}y_{t}+y_{t-5}y_{t})-(\frac{\alpha-\beta}{3})^{2}\lambda(E'_{(t-2)(t-1)},\vec{y})\nonumber\\
&=&\frac{\alpha-\beta}{3}(\frac{\alpha+2\beta}{3}+2\alpha)\frac{4\beta-\alpha}{3}-
\frac{(\alpha-\beta)^2}{9}\lambda(E'_{(t-2)(t-1)},\vec{y}).
\end{eqnarray}
Let $\vec{w}=(w_1, w_2,
\ldots, w_t)$ be given by $w_i=z_i$ for $i\neq t-4$, $i\neq t-3$ and
$w_{t-4}=z_{t-4}-\frac{\alpha-\beta}{3}=\frac{2\alpha+\beta}{3}$,
$w_{t-3}=z_{t-3}+\frac{\alpha-\beta}{3}=\frac{2\alpha+\beta}{3}$.
Clearly, $\vec{w}=(w_1,w_2, \ldots, w_t)$ is also a legal weighting
for $C_{3,m}$, and
\begin{eqnarray}\label{eq20e}
\lambda(C_{3,m},\vec{w})-\lambda(C_{3,m},\vec{z})
&=&\frac{\alpha-\beta}{3}[\lambda(E'_{t-3},\vec{z})-\lambda(E'_{t-4},\vec{z})]-(\frac{\alpha-\beta}{3})^{2}\lambda(E'_{(t-4)(t-3)},\vec{z})\nonumber\\
&=&\frac{\alpha-\beta}{3}(z_{t-4}-z_{t-3})\lambda(E'_{(t-4)(t-3)},\vec{z})-(\frac{\alpha-\beta}{3})^{2}\lambda(E'_{(t-4)(t-3)},\vec{z})\nonumber\\
&=&[\frac{\alpha-\beta}{3}(\alpha-\frac{\alpha+2\beta}{3})-(\frac{\alpha-\beta}{3})^{2}]\lambda(E'_{(t-4)(t-3)},\vec{z})\nonumber\\
&=&\frac{(\alpha-\beta)^2}{9}\lambda(E'_{(t-4)(t-3)},\vec{z}).
\end{eqnarray}
Let $\vec{u}=(u_1,
u_2, \ldots, u_t)$ be given by $u_i=w_i$ for $i\neq t-5$, $i\neq t-2$
and
$u_{t-5}=w_{t-5}-\frac{\alpha-\beta}{3}=\frac{2\alpha+\beta}{3}$,
$u_{t-2}=w_{t-2}+\frac{\alpha-\beta}{3}=\frac{2\alpha+\beta}{3}$.
Clearly, $\vec{u}=(u_1, u_2, \ldots, u_t)$ is also a legal weighting
for $C_{3,m}$, and
\begin{eqnarray}\label{eq21e}
\lambda(C_{3,m},\vec{u})-\lambda(C_{3,m},\vec{w})&=&\frac{\alpha-\beta}{3}[\lambda(E'_{t-2},\vec{w})-\lambda(E'_{t-5},\vec{w})]
-(\frac{\alpha-\beta}{3})^{2}\lambda(E'_{(t-5)(t-2)},\vec{w})\nonumber\\
&=&\frac{\alpha-\beta}{3}(w_{t-5}-w_{t-2})\lambda(E'_{(t-5)(t-2)},\vec{w})-(\frac{\alpha-\beta}{3})^{2}\lambda(E'_{(t-5)(t-2)},\vec{w})\nonumber\\
&=&[\frac{\alpha-\beta}{3}(\alpha-\frac{\alpha+2\beta}{3})-(\frac{\alpha-\beta}{3})^{2}]\lambda(E'_{(t-5)(t-2)},\vec{w})\nonumber\\
&=&\frac{(\alpha-\beta)^2}{9}\lambda(E'_{(t-5)(t-2)},\vec{w}).
\end{eqnarray}
Adding (\ref{g13x}), (\ref{eq18e}), (\ref{eq19e}), (\ref{eq20e}) and
(\ref{eq21e}), we have
\begin{eqnarray}\label{eq22e}
&&\lambda(C_{3,m},\vec{u})-\lambda(G,\vec{x})\nonumber\\
&=&\frac{(\alpha-\beta)^{2}}{9}[\lambda(E'_{(t-4)(t-3)},\vec{z})+\lambda(E'_{(t-5)(t-2)},\vec{w})- \lambda(E'_{(t-3)t},\vec{x})-\lambda(E'_{(t-2)(t-1)},\vec{y})]\nonumber\\
&& +\frac{\alpha-\beta}{3}(\beta^2+2\alpha\beta)+ \frac{\alpha-\beta}{3}(\frac{\alpha+2\beta}{3}+2\alpha)\frac{4\beta-\alpha}{3}+ 2(\beta-\alpha)\beta^2\nonumber \\
&=&\frac{(\alpha-\beta)^{2}}{9}(y_{t-2}+y_{t-1}+y_{t}+x_{t-3}+x_{t}+x_{t-1}-w_{t-2}-w_{t-5}-z_{t-3}-z_{t-4})+\frac{(\alpha-\beta)^2}{27}(37\beta-7\alpha)\nonumber \\
&=&\frac{(\alpha-\beta)^{2}}{9}(2\beta+\frac{4\beta-\alpha}{3}+3\beta-\frac{\alpha+2\beta}{3}-\alpha-\frac{\alpha+2\beta}{3}-\alpha)+\frac{(\alpha-\beta)^2}{27}(37\beta-7\alpha)\nonumber \\
&=&\frac{(\alpha-\beta)^2}{27}(52\beta-16\alpha)\nonumber \\
&\ge &0
\end{eqnarray}
since $\alpha\leq 2\beta$. Hence
$\lambda(C_{3,m})\geq\lambda(C_{3,m},\vec{u})\geq\lambda(G,\vec{x})=\lambda(G).$ \qed

\begin{lemma} \label{Lemma3} Let $G$ and $G'$ be  left-compressed $3$-graphs on vertex set $[t]$ with $m={t \choose 3}-a$ edges (where $7\le a\le t-2$) satisfying $|E(G) \triangle E(C_{3,m})|=|E(G')\triangle E(C_{3,m})|=6$ and the triples with the minimum colex ordering in $G^{c}$ and $G'^{c}$ are $(t-3)(t-2)(t-1)$ and $(t-5)(t-2)t$  respectively.  Then $ \lambda(G)\leq \lambda(G')\leq \lambda (C_{3,m})$.
\end{lemma}
{\em Proof.} By Lemma \ref{Lemma1}, $\lambda(G')\leq \lambda (C_{3,m})$. So it is sufficient to show $ \lambda(G)\leq \lambda(G')$.
Let $\vec{x}=(x_{1},x_{2},\ldots ,x_{t})$  be an optimal weighting for $G$ satisfying $x_1 \ge x_2 \ge \ldots \ge x_t \ge 0$.
We have
\begin{eqnarray}\label{eq101}
\lambda(G',\vec{x})-\lambda(G,\vec{x})=x_{t-3}x_{t-2}x_{t-1}-x_{t-5}x_{t-2}x_{t}.
\end{eqnarray}
Consider a new weighting ${\vec y}=(y_1, y_2, \ldots, y_t)$ given by $y_i=x_i$ for $i\neq t-5$, $i\neq t-3$ and $y_{t-5}=x_{t-5}-\delta$, $y_{t-3}=x_{t-3}+\delta$. Then
\begin{eqnarray}
\lambda(G',\vec{y})-\lambda(G',\vec{x})&=& \delta[\lambda(E'_{t-3},\vec{x})-\lambda(E'_{t-5},\vec{x})]-\delta^{2}\lambda(E'_{(t-5)(t-3)},\vec{y}) \nonumber\\
&=&\delta(x_{t-5}-x_{t-3})\lambda(E'_{(t-5)(t-3)},\vec{x})-\delta^{2}\lambda(E'_{(t-5)(t-3)},\vec{x}).
\end{eqnarray}
Let $\delta=\frac{x_{t-5}-x_{t-3}}{2}$. Clearly, ${\vec y}=(y_1, y_2, \ldots, y_t)$ is also a legal weighting for $G$. Also
\begin{eqnarray}\label{eq102}
\lambda(G',\vec{y})-\lambda(G',\vec{x})=\frac{(x_{t-5}-x_{t-3})^{2}}{4}\lambda(E'_{(t-5)(t-3)},\vec{x}).
\end{eqnarray}
Let ${\vec z}=(z_1, z_2, \ldots, z_t)$ be given by $z_i=y_i$ for $i\neq t-1$, $i\neq t$ and $z_{t-1}=y_{t-1}+\eta$, $z_{t}=y_{t}-\eta$. Then
\begin{eqnarray}
\lambda(G',\vec{z})-\lambda(G',\vec{y})&=& \eta[\lambda(E'_{t-1},\vec{y})-\lambda(E'_{t},\vec{y})]-\eta^{2}\lambda(E'_{(t-1)t},\vec{y})\nonumber\\
&=&\eta[y_{t-3}y_{t-2}+y_{t-4}y_{t-2}+y_{t-5}y_{t-2}-(y_{t-1}-y_{t})\lambda(E'_{(t-1)t},\vec{y})]-\eta^{2}\lambda(E'_{(t-1)t},\vec{y})\nonumber\\
&=&\eta[x_{t-3}x_{t-2}+x_{t-4}x_{t-2}+x_{t-5}x_{t-2}-(x_{t-1}-x_{t})\lambda(E_{(t-1)t},\vec{x})]-\eta^{2}\lambda(E_{(t-1)t},\vec{x})\nonumber\\
&=&\eta[x_{t-3}x_{t-2}+x_{t-4}x_{t-2}+x_{t-5}x_{t-2}-{\lambda(E_{(t-1)\backslash t},\vec{x})\over \lambda(E_{(t-1)t},\vec{x})} \lambda(E_{(t-1)t},\vec{x})]-\eta^{2}\lambda(E_{(t-1)t},\vec{x})\nonumber\\
&=&\eta(x_{t-3}x_{t-2}+x_{t-5}x_{t-2})-\eta^{2}\lambda(E_{(t-1)t},\vec{x}).
\end{eqnarray}
Let
$\eta=\frac{x_{t-3}x_{t-2}+x_{t-5}x_{t-2}}{2\lambda(E_{(t-1)t},\vec{x})}$.
By the condition of $|E(G) \triangle E(C_{3,m})|=6$ we have $\{1,2,3 \} \subseteq
E_{(t-1)t}$, so
$$\eta\le \frac{x_{t-2}}{3}.$$ Applying Remark \ref{r1}(b), we have
\begin{eqnarray}
x_{t-2}&=&x_{t}+{\lambda(E_{(t-2)\backslash t},\vec{x})\over
\lambda(E_{(t-2)t},\vec{x})}\nonumber\\
&\le &x_{t}+\frac{(x_{t-4}+\cdots + x_4)x_{t-1}}{1-x_t-x_{t-1}-x_{t-2}-x_{t-3}-x_{t-4}}\nonumber\\
&\leq&x_{t}+x_{t-1}\nonumber\\
&=&x_{t}+x_{t}+{\lambda(E_{(t-1)\backslash t},\vec{x})\over
\lambda(E_{(t-1)t},\vec{x})}\nonumber\\
&\leq&2x_{t}+{x_{t-4}x_{t-2}\over x_{1}+x_{2}+x_{3}}\nonumber\\
&\leq&2x_{t}+{x_{t-2}\over 3}.
\end{eqnarray}
So $\frac{x_{t-2}}{3}\leq x_{t}$. Recall that  $\eta\leq\frac{x_{t-2}}{3}$. Therefore, $\eta\leq x_{t}$. Hence ${\vec z}=(z_1,
z_2, \ldots, z_t)$ is also a legal weighting for $G$, and
\begin{eqnarray}\label{eq103}
\lambda(G',\vec{z})-\lambda(G',\vec{y})=\frac{(x_{t-5}+x_{t-3})^{2}x^{2}_{t-2}}{4\lambda(E_{(t-1)t},\vec{x})}.
\end{eqnarray}
By Remark \ref{r1}(b), we have
\begin{eqnarray}\label{eq104}
x_{t-5}=x_{t-3}+\frac{(x_{t-1}+x_{t})x_{t-2}}{\lambda(E_{(t-5)(t-3)},\vec{x})}.
\end{eqnarray}
Note that
  \begin{eqnarray}\label{eq105}\lambda(E_{(t-5)(t-3)},\vec{x})-\lambda(E_{(t-1)t},\vec{x})\geq (1-x_{t-5}-x_{t-3})-(1-x_{t-5}-x_{t-4}-x_{t-3}-x_{t-2}-x_{t-1}-x_{t})\geq 0.\end{eqnarray}
 %By Remark \ref{r1}(b), we have$$x_{t-5}=x_{t-2}+\frac{x_{t-3}x_{t-1}+x_{t-3}x_{t}+x_{t-1}x_{t}}{\lambda(E_{(t-5)(t-2)},\vec{x})}\leq x_{t-2}+x_{t-1}.$$Hence
 %\begin{eqnarray}\label{eq105}\lambda(E_{(t-5)(t-3)},\vec{x})\geq \lambda(E_{(t-1)t},\vec{x})\end{eqnarray}
Combing (\ref{eq101}), (\ref{eq102}),(\ref{eq103}),(\ref{eq104})and (\ref{eq105}), we have
\begin{eqnarray}
\lambda(G',\vec{z})-\lambda(G,\vec{x})&\geq&\frac{(x_{t-5}-x_{t-3})^{2}\lambda(E'_{(t-5)(t-3)},\vec{x})}{4} + \frac{(x_{t-5}+x_{t-3})^{2}x^{2}_{t-2}}{4\lambda(E_{(t-1)t},\vec{x})} \nonumber\\
&& -\frac{x_{t-2}^{2}(x_{t-1}+x_{t})x_{t}}{\lambda(E_{(t-5)(t-3)},\vec{x})}\nonumber\\
&=&\frac{x_{t-2}^{2}(x_{t-1}+x_{t})^2}{4\lambda(E_{(t-5)(t-3)},\vec{x})}+\frac{(x_{t-5}+x_{t-3})^{2}x^{2}_{t-2}}{4\lambda(E_{(t-1)t},\vec{x})}
-\frac{x_{t-2}^{2}(x_{t-1}+x_{t})x_{t}}{\lambda(E_{(t-5)(t-3)},\vec{x})}\nonumber\\&\ge&0.
\end{eqnarray}
Hence $\lambda(G)=\lambda(G,\vec{x})\leq \lambda(G',\vec{z})=\lambda(G').$  \qed

\begin{lemma} \label{Lemma4} Let $G$  be  the left-compressed $3$-graphs on vertex set $[t]$ with $m={t \choose 3}-6$ edges and $G^{c}=\{(t-2)(t-1)t,(t-3)(t-1)t,(t-4)(t-1)t,(t-3)(t-2)t,(t-4)(t-2)t,(t-3)(t-2)(t-1)\}$. Let $G^{'}=G\bigcup \{(t-4)(t-2)t\} \backslash\{(t-5)(t-1)t\}$. Then $ \lambda(G)\leq \lambda(G')\leq \lambda (C_{3,m})$.
\end{lemma}
{\em Proof.} By Lemma \ref{Lemma2}, $\lambda(G') \leq \lambda (C_{3,m})$. So it is sufficient to show that $\lambda(G)\leq \lambda(G')$. Let $\vec{x}=(x_{1},x_{2},\ldots ,x_{t})$  be an optimal weighting for $G$ satisfying $x_1 \ge x_2 \ge \ldots \ge x_t \ge 0$.
 Then
 \begin{eqnarray}\label{eq1050}
 \lambda(G',\vec{x})-\lambda(G,\vec{x})=x_{t-4}x_{t-2}x_{t}-x_{t-5}x_{t-1}x_{t}.
 \end{eqnarray}
 By Remark \ref{r1}(b), $x_{t-2}=x_{t-1}$. Consider a new weighting ${\vec y}=(y_1, y_2, \ldots, y_t)$ given by
 $y_i=x_i$ for $i\neq t-5$, $i\neq t-4$ and $y_{t-5}=x_{t-5}-\delta$, $y_{t-4}=y_{t-4}+\delta$. Then
\begin{eqnarray*}
\lambda(G',\vec{y})-\lambda(G',\vec{x})&=& \delta[\lambda(E'_{t-4},\vec{x})-\lambda(E'_{t-5},\vec{x})]-\delta^{2}\lambda(E'_{(t-5)(t-4)},\vec{x}) \nonumber\\
&=&\delta(x_{t-5}-x_{t-4})\lambda(E'_{(t-5)(t-4)},\vec{x})-\delta^{2}\lambda(E'_{(t-5)(t-4)},\vec{x}).
\end{eqnarray*}
Let $\delta=\frac{x_{t-5}-x_{t-4}}{2}$. Clearly, ${\vec y}=(y_1, y_2, \ldots, y_t)$ is also a legal weighting for $G$, and
\begin{eqnarray}\label{eq106}
\lambda(G',\vec{y})-\lambda(G',\vec{x})=\frac{(x_{t-5}-x_{t-4})^{2}}{4}\lambda(E'_{(t-5)(t-4)},\vec{x})=\frac{(x_{t-5}-x_{t-4})^{2}}{4}\lambda(E'_{(t-5)(t-4)},\vec{x}).
\end{eqnarray}
Let ${\vec z}=(z_1,z_2, \ldots, z_t)$ given by $z_m=y_m$ for $m\neq t-2$, $m\neq t-1$ and $z_{t-2}=y_{t-2}+\eta$, $z_{t-1}=y_{t-1}-\eta$. Then
\begin{eqnarray}
\lambda(G',\vec{z})-\lambda(G',\vec{y})&=& \eta[\lambda(E'_{t-2},\vec{y})-\lambda(E'_{t-1},\vec{y})]-\eta^{2}\lambda(E'_{(t-2)(t-1)},\vec{y}) \nonumber\\
&=&\eta(y_{t-4}y_{t}+y_{t-5}y_{t})-\eta^{2}\lambda(E'_{(t-2)(t-1)},\vec{y}).
\end{eqnarray}
Let
$\eta=\frac{y_{t-5}y_{t}+y_{t-4}y_{t}}{2\lambda(E'_{(t-2)(t-1)},\vec{y})}$.
Clearly, $\eta<y_{t}$. Hence, ${\vec z}=(z_1, z_2, \ldots, z_t)$ is
also a legal weighting for $G'$, and
\begin{eqnarray}\label{eq107}
\lambda(G',\vec{z})-\lambda(G',\vec{y})=\frac{(y_{t-5}+y_{t-4})^{2}y_{t}^{2}}{4\lambda(E'_{(t-2)(t-1)},\vec{y})}=\frac{(x_{t-5}+x_{t-4})^{2}x_{t}^{2}}{4\lambda(E_{(t-2)(t-1)},\vec{x})}.
\end{eqnarray}
By Remark \ref{r1}(b), we have
\begin{eqnarray}\label{eq108}
x_{t-5}=x_{t-4}+\frac{2x_{t-1}x_{t}}{\lambda(E_{(t-5)(t-4)},\vec{x})}\le x_{t-4}+x_{t},
\end{eqnarray}
 \begin{eqnarray}\label{eq110}
 x_{t-4}=x_{t-1}+\frac{\lambda(E_{(t-4)\backslash (t-1)},\vec{x})}{\lambda(E_{(t-4)(t-1)},\vec{x})}= x_{t-1}+\frac{x_{t-3}x_{t}+x_{t-3}x_{t-2}}{1-x_{t-4}-x_{t-1}-x_{t}},\end{eqnarray}
 \begin{eqnarray}\label{eq111}x_{t-5}=x_{t-1}+\frac{\lambda(E_{(t-5)\backslash (t-1)},\vec{x})}{\lambda(E_{(t-5)(t-1)},\vec{x})}= x_{t-1}+\frac{(x_{t-4}+x_{t-3}+x_{t-2})x_{t}+x_{t-3}x_{t-2}}{1-x_{t-5}-x_{t-1}}.\end{eqnarray}
 Combing (\ref{eq108}), (\ref{eq110}) and (\ref{eq111}) we have
\begin{eqnarray*}x_{t-4}+x_{t-5}&\le& 2x_{t-1}+\frac{(x_{t-4}+x_{t-3}+x_{t-2}+x_{t-3})x_{t}+2x_{t-3}x_{t-2}}{1-x_{t-4}-x_{t-1}-x_{t}}\\
 &\le&2x_{t-1}+x_{t}+x_{t-3}.\end{eqnarray*}
 Hence,
\begin{eqnarray}\label{eq109}
\lambda(E_{(t-5)(t-4)},\vec{x})-\lambda(E_{(t-2)(t-1)},\vec{x})= x_{t-3}+2x_{t-1}+x_{t}-x_{t-4}-x_{t-5}\geq 0.
\end{eqnarray}
Combing (\ref{eq1050}), (\ref{eq106}), (\ref{eq107}), (\ref{eq108}) and (\ref{eq109}),  we have
\begin{eqnarray*}
\lambda(G',\vec{z})-\lambda(G,\vec{x})&=\frac{x_{t-1}^{2}x_{t}^{2}}{\lambda(E_{(t-5)(t-4)},\vec{x})}
+\frac{(x_{t-5}+x_{t-4})^{2}x^{2}_{t}}{4\lambda(E_{(t-2)(t-1)},\vec{x})}
-\frac{2x_{t-1}^{2}x_{t}^{2}}{\lambda(E_{(t-5)(t-4)},\vec{x})}\geq 0.
\end{eqnarray*}
Hence $\lambda(G')\geq\lambda(G',\vec{z})\geq\lambda(G,\vec{x})=\lambda(G).$ \qed

Now we are ready to show Theorem \ref{mainresult}.

\bigskip
\noindent{\em Proof of Theorem \ref{mainresult}.}  Let $G$ be a left-compressed extremal
$3$-graph with $m$ edges satisfying $|E(G) \Delta E(C_{3,m})| \le 6$.  Let $t$ be a positive integer such that ${t-1 \choose 3} \le m < {t \choose 3}$. By Remark \ref{reducetot},  we can assume that $G=([t], E)$. If  $m<{t-1 \choose 3}+{t-2 \choose 2}$, let $G'$ be obtained by  adding the first  ${t-1 \choose 3}+{t-2 \choose 2}-m$ triples  in colex ordering in $E(G^c)$ to $E(G)$, then $\lambda(E(G'))\ge\lambda(E(G))$ and $|E(G') \Delta E(C_{3,m'})| \le 6$ for $m'={t-1 \choose 3}+{t-2 \choose 2}$.  In view of  Lemma \ref{LemmaTal7},
we may assume  that $m\ge {t-1 \choose 3}+{t-2 \choose 2}$. Then  $C_{3,m}\supseteq
[t-1]^{(3)}\bigcup \{1,2,\cdots,t-2,t\}^{(3)}$.  Since $|E \Delta E(C_{3,m})| \le 6$, there are at most 3 triples in $\{1,2,\cdots,t-2,t\}^{(3)}\bigcup [t-1]^{(3)}-E$. Since $G$ is left-compressed, in view of Figure 1, there are only the following possible cases for $\{1,2,\cdots,t-2,t\}^{(3)}\bigcup [t-1]^{(3)}-E$:

Case 1. $\{1,2,\cdots,t-2,t\}^{(3)}\bigcup [t-1]^{(3)}-E=\{(t-3)(t-2)t\}$;

 Case 2. $\{1,2,\cdots,t-2,t\}^{(3)}\bigcup [t-1]^{(3)}-E=\{(t-3)(t-2)t,(t-4)(t-2)t \}$;

Case 3. $\{1,2,\cdots,t-2,t\}^{(3)}\bigcup [t-1]^{(3)}-E=\{(t-3)(t-2)t,(t-3)(t-2)(t-1)\}$;

 Case 4. $\{1,2,\cdots,t-2,t\}^{(3)}\bigcup [t-1]^{(3)}-E=\{(t-3)(t-2)t,(t-4)(t-2)t, (t-5)(t-2)t \}$;

Case 5. $\{1,2,\cdots,t-2,t\}^{(3)}\bigcup [t-1]^{(3)}-E=\{(t-3)(t-2)t,(t-4)(t-2)t, (t-4)(t-3)t \}$;

Case 6. $\{1,2,\cdots,t-2,t\}^{(3)}\bigcup [t-1]^{(3)}-E=\{(t-3)(t-2)t, (t-4)(t-2)t,(t-3)(t-2)(t-1)\}$.

If Cases 1, 2, 4, 5 happen, then by Lemma \ref{Lemma1}, $\lambda(G)\leq\lambda(C_{3,m})$. If Case 3 and $a \ge 5$ happen, then by Lemma \ref{Lemma2}, $\lambda(G)\leq\lambda(C_{3,m})$. If Case 3 and $a < 5$ happen, then $a=4$. By Lemma \ref{Lemma2+}, $\lambda(G)\leq\lambda(C_{3,m})$. If Case 6 and $a\ge 7$ happen, then by Lemma \ref{Lemma3}, $\lambda(G)\leq\lambda(C_{3,m})$. If Case 6 and $a< 7$ happen, then $a=6$.  By Lemma \ref{Lemma4}, $\lambda(G)\leq\lambda(C_{3,m})$. The proof of Theorem \ref{mainresult} is completed.
\qed

\bigskip
\noindent{\em Proof of Corollary \ref{coro1}.} We can assume $G$ is a left-compressed $3$-graph on $[t]$ by Remark \ref{reducetot}.  The range of $m$ gurantees that
$|E(G) \Delta E(C_{3,m})|\leq 6 $. Therefore $\lambda(G)\leq \lambda (C_{3,m})$ by Theorem \ref{mainresult}. \qed

The following result is also implied by Theorem \ref{mainresult}.
\begin{coro} Let $m$, $t$ and $b$ be positive integers satisfying
$${t-1 \choose 3}+ {t-2 \choose 2}+b< m \leq{t \choose 3}.$$  Let $G=(V,E)$ be a left-compressed 3-graph on the vertex set [t] with $m$ edges satisfying $|E_{(t-1)t}|\le b+3$. Then $\lambda (G)\leq \lambda (C_{3,m})$.
\end{coro}
 {\em Proof.} Since $|E_{(t-1)t}|\le b+3$, we have  $|E\Delta E(C_{3,m})|\leq 6$. So $\lambda (G)\leq \lambda (C_{3,m})$ by Theorem \ref{mainresult}. \qed

\bigskip
{\bf Acknowledgments.} This research is partially  supported by National Natural Science Foundation of China (No. 11271116)

\end{document}